\documentclass[a4paper,11pt]{article}

\usepackage[utf8]{inputenc}
\usepackage[english,activeacute]{babel}
\usepackage{amsmath,amsthm,amsfonts,amssymb}
\usepackage{mathpartir}
\usepackage{ stmaryrd }
\usepackage{graphics}
\usepackage{enumerate}
\usepackage{mathrsfs}
\usepackage{rotating}
\usepackage{hyperref}
\usepackage{xcolor}
\hypersetup{
    colorlinks,
    linkcolor={red!50!black},
    citecolor={blue!50!black},
    urlcolor={blue!80!black}
}
\input{xy}
\usepackage[all]{xy}
\usepackage{tikz}
\usepackage{tikz-cd}

\newtheorem{thm}{Theorem}[section]
\newtheorem*{thm*}{Theorem}

\newtheorem{lemma}[thm]{Lemma}

\theoremstyle{plain}

\theoremstyle{remark}

\input{diagxy}
\xyoption{curve}

\newbox\anglebox % large pullback angle
\setbox\anglebox=\hbox{\xy \POS(75,0)\ar@{-} (0,0) \ar@{-} (75,75)\endxy}
 
\newbox\angleboxr % reverse large pullback angle
\setbox\angleboxr=\hbox{\xy \POS(0,0)\ar@{-} (0,75) \ar@{-} (75,0)\endxy}
 
\newbox\sanglebox % small pullback angle
\setbox\sanglebox=\hbox{\xy \POS(50,0)\ar@{-} (0,0) \ar@{-} (50,50)\endxy}
 
\newbox\sangleboxr % small reverse pullback angle
\setbox\sangleboxr=\hbox{\xy \POS(0,0)\ar@{-} (0,50) \ar@{-} (50,0)\endxy}
 
\newbox\sangleboxf % small flipped pullback angle
\setbox\sangleboxf=\hbox{\xy \POS(50,50)\ar@{-} (50,0) \ar@{-} (0,50)\endxy}
 
\newbox\angleboxf % flipped pullback angle
\setbox\angleboxf=\hbox{\xy \POS(75,75)\ar@{-} (75,0) \ar@{-} (0,75)\endxy}
 
\newbox\sangleboxfr % small flipped reverse pullback angle
\setbox\sangleboxfr=\hbox{\xy \POS(0,50)\ar@{-} (50,50) \ar@{-} (0,0)\endxy}
 
\newbox\angleboxfr % small flipped reverse pullback angle
\setbox\angleboxfr=\hbox{\xy \POS(0,75)\ar@{-} (75,75) \ar@{-} (0,0)\endxy}

\newcommand{\Sets}{\ensuremath{\mathbf{Set}}}

\title{A complete classification of categoricity spectra of accessible categories with directed colimits}
\author{Christian Esp\'indola}

\begin{document}
\date{}
\maketitle

\begin{abstract}
We provide a complete classification of all the possible categoricity spectra, in terms of internal size, that can appear in a large accessible category with directed colimits, assuming the Singular Cardinal Hypothesis ($SCH$), and providing as well explicit threshold cardinals for eventual categoricity. This includes as a particular case the first complete classification of categoricity spectra of abstract elementary classes (AEC's) entirely in $ZFC$. More specifically, we have:

\begin{thm*} Let $\mathcal{K}$ be a large $\kappa$-accessible category with directed colimits. Assume the Singular Cardinal Hypothesis $SCH$ (only if the restriction to monomorphisms is not an AEC). Then the categoricity spectrum $\mathcal{C}at(\mathcal{K})=\{\lambda\geq \kappa: \mathcal{K} \text{ is $\lambda$-categorical}\}$ is one of the following:

\begin{enumerate}
\item $\mathcal{C}at(\mathcal{K})=\emptyset$.
\item $\mathcal{C}at(\mathcal{K})=[\alpha, \beta]$ for some $\alpha, \beta \in [\kappa, \beth_{\omega}(\kappa))$.
\item $\mathcal{C}at(\mathcal{K})=[\chi, \infty)$ for some $\chi \in [\kappa, \beth_{(2^{\kappa})^+})$.
\end{enumerate}

\end{thm*}

This solves in particular Shelah categoricity conjecture for AEC's. There are examples of each of the three cases of the classification, showing that they indeed occur. 
\end{abstract}

\section{Introduction}

This short paper is a sequel to the work of the author in \cite{espindolas} in which a generalization of Shelah's eventual categoricity conjecture (Conjecture 4.2 in the introduction of \cite{shelah}) is proven in the more general context of accessible categories with directed colimits. When all morphisms are monomorphisms in such categories of models, an analogous form of Shelah's presentation theorem exhibits them as a projective class of an infinite quantifier logic, for which even the Hanf number for model existence has no known explicit bound in $ZFC$ (the only known bound is a strongly compact cardinal, and in fact in some models of $ZFC$ the Hanf number for $\mathcal{L}_{\omega_1, \omega_1}$ exceeds the first measurable cardinal).

In the special case of those accessible categories which have directed colimits, however, we will prove that one can find $ZFC$ bounds for the Hanf number of model existence. Grossberg has emphasized the importance of also having explicit threshold cardinals for the eventual categoricity phenomenon, and we now intend to use the same setup and results of \cite{espindolas}, together with Morley's method, to provide the provably best possible explicit thresholds. Shelah categoricity conjecture asks to prove, in $ZFC$, that the threshold for eventual categoricity in an AEC $\mathcal{K}$ is $\beth_{(2^{LS(\mathcal{K})})^+}$ (see Conjecture 4.3 b) in the introduction of \cite{shelah}). We will prove this conjecture here. An example from Shelah mentioned in \cite{sv3} shows that this threshold is best possible.

   Assuming $SCH$, we are also going to provide a proof of a direct generalization of Shelah categoricity conjecture to the more general context of accessible categories with directed colimits. If $\mathcal{K}$ is such a category, we show that such that if $\mathcal{K}$ is categorical in some in some $\lambda \geq \beth_{(2^{LS(\mathcal{K})})^+}$ (i.e., it has only one object of some high enough internal size up to isomorphism), then $\mathcal{K}$ is $\lambda'$-categorical for every $\lambda' \geq \beth_{(2^{LS(\mathcal{K})})^+}$. When considering cardinalities of models of infinitary theories $\mathbb{T}$ of $\mathcal{L}_{\kappa, \theta}$ that axiomatize $\mathcal{K}$, the result implies, under $SCH$, the following infinitary version of Morley's categoricity theorem:

\begin{thm}\label{morley}
(Morley's categoricity theorem for infinitary theories) Let $\phi$ be a $\mathcal{L}_{\kappa, \theta}$ sentence whose category of models and $\mathcal{L}_{\kappa, \theta}$-elementary embeddings has directed colimits. Let $S$ be the class of cardinals $\lambda$ which are of cofinality at least $\theta$ but are not successors of cardinals of cofinality less than $\theta$. Assume the weakening of the Singular Cardinal Hypothesis $SCH_{\theta, \geq 2^{<\theta}}$. Then, if $\phi$ is $\lambda$-categorical for some $\lambda \geq \beth_{(2^{\kappa})^+}$ in $S$, then $\phi$ is $\lambda'$-categorical for every $\lambda' \geq \beth_{(2^{\kappa})^+}$ in $S$. Moreover:
\begin{enumerate}
\item if the directed colimits are concrete, we can spare the assumption $SCH_{\theta, \geq 2^{<\theta}}$ and take $S$ as the class of all cardinals.
\item if $\phi$ is compact and the morphisms of our category are $\mathcal{L}_{\omega, \omega}$-elementary embeddings, we can replace $\beth_{(2^{\kappa})^+}$ with $\kappa$.
\end{enumerate}
\end{thm}

Here $SCH_{\theta, \geq 2^{<\theta}}$ is defined as ``for all $\mu \geq 2^{<\theta}$ there is a set of cardinals $\lambda_i \leq \mu$ unbounded below $\mu$ such that, for each $i$, $\nu^{<\theta} \leq \lambda_i$ for all $\nu<\lambda_i$", see Remark 2.3 of \cite{rlv}.  Also, we know from \cite{espindolas} that there are examples showing that the exceptions in the class $S$ are needed.

The case $\theta=\omega$ in Theorem \ref{morley} is Shelah categoricity conjecture for $\mathcal{L}_{\kappa, \omega}$, since in this case $SCH_{\omega, \geq 2^{<\omega}}$ is provable in $ZFC$. When the directed colimits are concrete, since we restrict to monomorphisms, the result is Shelah categoricity conjecture for AEC's, since $SCH_{\theta, \geq 2^{<\theta}}$ can also be removed by the methods of \cite{espindolas}. The case when $\phi$ is compact (i.e., when it has the property that $\phi$ is consistent with an arbitrary set of first-order finitary formulas if and only if it is consistent with each of its finite subsets, see \cite{keisler}) is precisely a proper generalization of Morley's categoricity theorem (for countable first-order theories) and of Shelah's categoricity theorem (for uncountable first-order theories). Compact sentences of infinitary logic forms a much wider class than these two particular cases, since they include (but are not limited to) all those conjunctive sentences (i.e., sentences where only conjunctions are infinitary but disjunctions are finitary). Thus, Theorem \ref{morley} is a vast generalization of those conjectures and results to the realm of infinite quantifier theories and provides new proofs of the known theorems for finitary first-order theories. 

  As it turns out, the existence of directed colimits is what allows for a smooth classification theory. The main tool for this will be Theorem \ref{tameness}, which at the same time extends work of Shelah for $\mathcal{L}_{\omega_1, \omega}$ showing, under the Weak Generalized Continuum Hypothesis ($WGCH$) that categoricity in the first $\omega$ cardinals implies categoricity everywhere (see \cite{shelah}). We remove here the set-theoretic hypothesis and generalize this result to AEC's, at the price of asking for categoricity in the first $\beth_{\omega}$ cardinals. By the example of Shelah and Villaveces in \cite{shvi}, this seems to be close to optimal, since they showed that categoricity can fail above $\beth_n(\lambda)$ for any $n \in \omega$ while holding at the first $n$ cardinals above the Löwenheim-Skolem number $\lambda$ (though it is open whether categoricity holds up to $\beth_n(\lambda)$ or the gap could be reduced further). Theorem \ref{tameness} also uses higher dimensional amalgamation properties, which are shown to be a consequence of categoricity by means of a simple categorical proof, thereby simplifying the methods of \cite{shva}.
  
  Finally, we state the classification of categoricity spectra in AEC's, in $ZFC$, and assuming $SCH$ also in accessible categories with directed colimits (the set theoretic assumption is needed to guarantee that the existence spectrum contains an end tail of cardinals). This uses Lemma \ref{dc}, some weaker versions of which in the context of AEC's have appeared in the literature. We give here a categorical proof based, among other things, on a form of Lawvere's duality for algebraic theories. The proof extends the result to $\mu$-AEC's with directed colimits, and greatly simplifies the arguments given for AEC's, in such a way that it can be applied to deduce eventual categoricity without needing to use amalgamation, using instead an observation on the double negation topology.
  
  As a word of warning, we emphasize that all the methods, results and notation from the authors' previous paper \cite{espindolas} are assumed here throughout, so the reader is advised to go through that paper first before continuing with this sequel.

\section{Saturation and stability}

We start by stating the following lemma of independent interest:

\begin{lemma}\label{dc}
Let $\mathcal{K}$ be a $\mu$-AEC with directed colimits and amalgamation that is $\rho$-stable for each $\rho<\kappa$. Then $\kappa$-saturated models are closed under directed colimits.
\end{lemma}

\begin{proof}
Consider the topos $\Sets[\mathbb{T}_{\kappa}^B]_{\lambda}/[M, -] \cong \Sets^{\mathcal{K}_{\geq \kappa, < \lambda}^B}/[M, -] \cong \Sets^{M/\mathcal{K}_{\geq \kappa, < \lambda}^B}$, (where $\mathcal{K}_{\geq \kappa, < \lambda}^B$ consists of the models in $\mathcal{K}_{\geq \kappa, < \lambda}$ and all its $\kappa$-Boolean homomorphisms, and where $\mathbb{T}_{\kappa}^B$ is $\mathbb{T}_{\kappa}$ plus all those instances of excluded middle for $\kappa$-coherent formulas). We have a stable surjection $\Sets^{M/\mathcal{K}_{\geq \kappa, < \lambda}^B} \twoheadrightarrow \Sets^{M/\mathcal{K}_{\geq \kappa, < \lambda}}$; this can be seen by considering first the stable surjection $\Sets^{\mathcal{K}_{\geq \kappa, < \lambda}^B} \cong  \Sets[\mathbb{T}_{\kappa}^B]_{\lambda} \twoheadrightarrow \Sets[\mathbb{T}_{\kappa}]_{\lambda} \cong \Sets^{\mathcal{K}_{\geq \kappa, < \lambda}}$. Then we consider the pullback functor to the slice $\Sets^{\mathcal{K}_{\geq \kappa, < \lambda}} \to \Sets[\mathbb{T}_{\kappa}]_{\lambda}/[M, -] \cong \Sets^{M/\mathcal{K}_{\geq \kappa, < \lambda}}$, which is a geometric morphism along whose direct image we take the following (pseudo-)pullback:

\begin{center}
\begin{tikzcd}
\Sets^{\mathcal{K}_{\geq \kappa, < \lambda}^B} \arrow[rr, two heads]                                   &  & \Sets^{\mathcal{K}_{\geq \kappa, < \lambda}}                                   \\
                                                                                         &  &                                                                  \\
{\Sets[\mathbb{T}_{\kappa}^B]_{\lambda}/[M_i, -]} \arrow[uu] \arrow[rr, two heads] &  & {\Sets[\mathbb{T}_{\kappa}]_{\lambda}/[M_i, -]} \arrow[uu] \\
                                                                                         &  &                                                                  \\
{\Sets[\mathbb{T}_{\kappa}^B]_{\lambda}/[M, -]} \arrow[rr, two heads] \arrow[uu]   &  & {\Sets[\mathbb{T}_{\kappa}]_{\lambda}/[M, -]} \arrow[uu]  
\end{tikzcd}
\end{center}

The (pseudo-)pullback is precisely $\Sets[\mathbb{T}_{\kappa}^B]_{\lambda}/[M, -] \cong \Sets^{M/\mathcal{K}_{\geq \kappa, < \lambda}^B}$, as can be verified using the universal property of the slice. More generally, if we have now a sequence of embeddings $M_0 \to M_1 \to ... \to M$ with directed colimit $M$, then $\Sets^{M/\mathcal{K}_{\geq \kappa, < \lambda}}$ will be the limit of the chain formed by the $\Sets^{M_i/\mathcal{K}_{\geq \kappa, < \lambda}}$ and induced by those embeddings. Since pullbacks preserve limits, this implies that $\Sets^{M/\mathcal{K}_{\geq \kappa, < \lambda}^B}$ will be the limit of the chain formed by the $\Sets^{M_i/\mathcal{K}_{\geq \kappa, < \lambda}^B}$; in particular (considering functors from the presheaves to $\Sets$ preserving limits and colimits), the Cauchy completion of the slice $M/\mathcal{K}_{\geq \kappa, < \lambda}^B$ is the (pseudo-)limit in $\mathcal{C}at$ of the Cauchy completion of the slices $M_i/\mathcal{K}_{\geq \kappa, < \lambda}^B$ (note, in turn, that the Cauchy completions of the slices are equivalent to the slices of the Cauchy completion $\overline{\mathcal{K}_{\geq \kappa, < \lambda}^B}$).

  Note that for any $\kappa$-small model $P$ and $\kappa$-geometric theory $\mathbb{S}_{\kappa}$ of models of size at least $\kappa$ containing $P$, the pullback, in the $2$-category of $\kappa^+$-toposes and $\kappa^+$-geometric morphisms, of the double negation subtopos $\mathcal{S}$ of $\Sets[\mathbb{S}_{\kappa}]_{\kappa^+}$ along the surjection $s: \Sets[\mathbb{S}^B_{\kappa}]_{\kappa^+} \twoheadrightarrow \Sets[\mathbb{S}_{\kappa}]_{\kappa^+}$ must be $\mathcal{S}$ itself, so $\mathcal{S}$ embeds into $\Sets[\mathbb{S}^B_{\kappa}]_{\kappa^+}$. Moreover, the embedding is dense: take a nonzero subterminal object $A$ in $\Sets[\mathbb{S}^B_{\kappa}]_{\kappa}$; then it is nonzero in some model $M$ which, by density, must embed into the model of size $\kappa$. Then, if $\mathcal{M}$ is the category of models of $\mathbb{S}$, the $\kappa^+$-coherent sentence $\psi=\lim ev_{\phi_i}$ in $\Sets^{\mathcal{M}^B_{\kappa}}$, where $\lim ev_{\phi_i} \cong [M, -]$ in $\Sets^{\mathcal{M}_{\kappa}}$, is non zero in $\Sets^{\mathcal{M}_{\kappa}^{sat}}$ and it implies $A$, which is thus also non zero (and thus nonzero in $\mathcal{S}$, as we claimed). Note also that, as a side consequence of this proof, double negation commutes with $\kappa^+$-small conjunctions in $\Sets^{\mathcal{M}_{\kappa}^B}$.
  
  Assume now that all $M_i$ are $\kappa$-saturated. Without loss of generality we can also assume that $\kappa=\delta^+$ is a successor, since for limit $\kappa$ the saturated model is a directed colimit of smaller saturated models. Let us now prove that $M$ must be $\kappa$-closed (whence also $\kappa$-saturated). So consider an embedding $f: M \to N$; since $\mathcal{K}$ is $\rho$-stable with respect to Galois types over some $\delta$-saturated submodel $P$, it is $\rho$-stable with respect to $\kappa$-Boolean types of the same kind, so that an application of the omitting types theorem from \cite{espindolas} to the $\kappa$-Boolean theory of models of size at least $\kappa$ containing $P$, $\mathbb{S}^B_{\kappa}$, shows that all subobject lattices of $(\mathbf{x}, \top)$ in $\Sets[\mathbb{S}^B_{\kappa}]_{\kappa}$, for $\mathbf{x}$ a nonempty finite tuple, are atomic and thus Boolean. On the other hand, note that $\Sets[\mathbb{S}_{\kappa}]_{\kappa}$ is two-valued since it is a subtopos of $\Sets[\mathbb{S}_{\delta}]_{\kappa}$ and that this latter is equivalent to the slice $\Sets^{\mathcal{K}_{\delta}}/\phi_P$, where $\phi_P$ is the $\kappa$-geometric existential sentence corresponding to the diagram of $P$, which is two-valued since $\phi_P$ is an atom, as every model from $\mathcal{K}_{\delta}$ embeds in $P$. Therefore, the mentioned subobject lattices in $\Sets[\mathbb{S}^B_{\kappa}]_{\kappa}$ coincide with those in $\Sets[\mathbb{S}_{\kappa}]_{\kappa}$, and this entails, in particular, that $\Sets[\mathbb{S}^B_{\kappa}]_{\kappa}$ is two-valued. This readily implies that the colimit coprojections as well as their composition with $f$ are $\kappa$-Boolean, and so we have $\kappa$-Boolean embeddings $M_i \to N$, which induce a cone between the slices $N/\overline{\mathcal{K}_{\geq \kappa, < \lambda}^B}$ and $M_i/\overline{\mathcal{K}_{\geq \kappa, < \lambda}^B}$. By the universal property of the limit, there is an induced functor $N/\overline{\mathcal{K}_{\geq \kappa, < \lambda}^B} \to M/\overline{\mathcal{K}_{\geq \kappa, < \lambda}^B}$, which provides a natural transformation $[N, -] \to [M, -]$. By Yoneda, this must correspond to a morphism $M \to N$ in $\mathcal{K}_{\geq \kappa, < \lambda}^B$, and since this must be $f$, it follows that $f$ is $\kappa$-closed, as we wanted to show.
\end{proof}

\section{Categoricity and tameness}

We start by showing that in any $\mu$-AEC with amalgamation and no maximal models, categoricity in a high enough cardinal implies eventual tameness. We consider the same setup of section 8 in \cite{espindolas}, which we reproduce for the sake of convenience. Given a $\mu$-AEC $\mathcal{K}$ with Löwenheim-Skolem number $\kappa$ and $\mu \leq \kappa^+$, following Baldwin-Boney-Vasey, we add a $\kappa^+$-small arity predicate $P$ whose interpretation in a model $M$ consists of the image of the underlying structure of a model $N$ of size $\kappa$ embedded in $M$ through a morphism in the $\mu$-AEC. This particular expansion, which gives rise to an isomorphic AEC, has the property that morphisms coincide with substructure embeddings. Moreover, its models of size at least $\kappa$ can be axiomatized as follows, extending further the language with the symbol $\subseteq$:

$$\top \vdash_{\mathbf{x}} \exists \mathbf{y}\left(\bigvee_{M_0 \in S}\psi_{M_0}(\mathbf{y}) \wedge \mathbf{x} \subseteq \mathbf{y} \wedge P(\mathbf{y})\right)$$

$$\top \vdash_{\mathbf{x}\mathbf{y}} (\mathbf{x} \subseteq \mathbf{y} \wedge P(\mathbf{x}) \wedge P(\mathbf{y})) \to \bigvee_{(M_0, M_1) \in T}\psi_{(M_0, M_1)}(\mathbf{x}, \mathbf{y})$$

$$\top \vdash_{\mathbf{x}\mathbf{y}} (P(\mathbf{y}) \wedge \psi_{(M_0, M_1)}(\mathbf{x}, \mathbf{y})) \to P(\mathbf{x})$$

$$\top \vdash_{\mathbf{x}\mathbf{y}} \mathbf{x} \subseteq \mathbf{y} \leftrightarrow \bigwedge_{i \in I}\bigvee_{j \in J}x_i=y_j$$
\noindent
Here $S$ is a skeleton of the subcategory of models of size $\kappa$, $T$ is the set of pairs $(M_0, M_1)$
with a morphism in the $\mu$-AEC and $M_0, M_1 \in S$, while $\psi_{M_0}, \psi_{M_0, M_1}$ are conjunctions of
atomic and negated atomic formulas of the extended language such that $\psi_{M_0}(\mathbf{z})$ holds if and only if $\mathbf{z}$ is isomorphic to $M_0$, and $\psi_{M_0, M_1}(\mathbf{z}, \mathbf{w})$ holds if and only if $(\mathbf{z}, \mathbf{w})$ is isomorphic to $(M_0, M_1)$. 

  Assuming now categoricity at $\kappa$, we can get an axiomatization of an isomorphic $\mu$-AEC which can be entirely rewritten through sequents in the $(2^{\kappa})^+$-$Reg_{\neg}$ fragment. This is an intuitionistic fragment of first-order logic which contains no disjunctions, obtained from the $(2^{\kappa})^+$-regular fragment by adding $\bot$, together with the axioms $\bot \vdash_{\mathbf{x}} \phi$ and the axioms for $\neg$ that make it into a negation operator. Indeed, in the first sequent above the disjunction reduces to a single disjunct since we have categoricity at $\kappa$, while the last three sequents above have the general form of universal sentences $\forall \mathbf{z} \bigvee_{i \in I} \bigwedge_{j \in J} \psi_{ij}$, and each such sentence is equivalent to the set of sequents $\{\exists \mathbf{z} \bigwedge_{i \in I} \neg \psi_{if(i)} \vdash \bot\}_{f \in J^I}$.

  The $(2^{\kappa})^+$-$Reg_{\neg}$ fragment contains the $(2^{\kappa})^+$-$Reg_{\bot}$ subfragment, not containing the symbol $\neg$. The syntactic category $\mathcal{C}$ of any $(2^{\kappa})^+$-$Reg_{\neg}$ theory can be studied through the category $\mathcal{K}_{\geq (2^{\kappa})^+}^r$ of its $(2^{\kappa})^+$-$Reg_{\bot}$ models (models of the $(2^{\kappa})^+$-$Reg_{\bot}$ internal theory of $\mathcal{C}$, also known as the $(2^{\kappa})^+$-$Reg_{\bot}$ Morleyization of the $(2^{\kappa})^+$-$Reg_{\neg}$ theory). These latter are in particular $(2^{\kappa})^+$-regular models for the extended signature in which there is an extra propositional symbol $\bot$ and one predicate symbol $S$ for each negated atomic formula $\neg R$ and where the axioms of the theory contain all axioms obtained from formally replacing $\neg R$ by $S$ in each $(2^{\kappa})^+$-$Reg_{\neg}$ axiom and, additionally, all those axioms of the form $\bot \vdash_{\mathbf{x}} \phi$ and $R \wedge S \vdash_{\mathbf{x}} \bot$.

   If $(\mathcal{C})_{\lambda^+}^r$ is the syntactic category of the $\lambda^+$-$Reg_{\bot}$ theory with the same axioms as the $(2^{\kappa})^+$-$Reg_{\bot}$ theory of $\mathcal{C}$, then its $\lambda^+$-classifying topos $\mathcal{S}h((\mathcal{C})_{\lambda^+}^r, \tau)$ (where $\tau$ is the $\lambda^+$-$Reg_{\bot}$ coverage) will be precisely equivalent to the presheaf topos $\Sets^{\mathcal{K}_{\geq (2^{\kappa})^+, \leq \lambda}^r}$, as can be seen as a special case of Theorem 4.1 from \cite{espindolas} when $\lambda$ is big enough. In particular, the embedding $(\mathcal{C})_{\lambda^+}^r \to \Sets^{\mathcal{K}_{\geq (2^{\kappa})^+, \leq \lambda}^r}$ will preserve $\neg$ since it can be identified with Yoneda embedding, which preserves any right adjoint to pullback functors that might exist, see \cite{bj}).

  Using the compactness of $(2^{\kappa})^+$-$Reg_{\bot}$ logic, it is also easy to verify that the canonical functor $F: \mathcal{C} \to (\mathcal{C})_{\lambda^+}^r$ also preserves $\neg$. For if given a $\lambda^+$-regular sentence $\exists \mathbf{x} \bigwedge_{i<\lambda}\phi_i$ we have $\exists \mathbf{y} \bigwedge_{i<\lambda}\phi_i \wedge R \vdash_{\mathbf{x}} \bot$ in $\lambda^+$-$Reg_{\bot}$ logic, there must be a $(2^{\kappa})^+$-regular sentence $\exists_{i \in T} y_i \bigwedge_{i \in T}\phi_i$, for some subset $T \subset \lambda$ of size at most $2^{\kappa}$, such that $\exists_{i \in T} y_i \bigwedge_{i \in T}\phi_i \wedge R \vdash_{\mathbf{x}} \bot$ in $(2^{\kappa})^+$-$Reg_{\bot}$ logic, from which our result follows.

  It follows, in fact, that the evaluation functor $ev: \mathcal{C} \to \Sets^{\mathcal{K}_{\geq (2^{\kappa})^+, \leq \lambda}^r}$, the composite of Yoneda embedding with $F$, preserves $\neg$,\footnote{It is also possible to give a direct proof of this fact, using the compactness of $(2^{\kappa})^+$-$Reg_{\bot}$ logic, with the same arguments as in the proof of Joyal's theorem, according to which $ev: \mathcal{C} \to \Sets^{Mod(\mathcal{C})}$ preserves universal quantification when $Mod(\mathcal{C})$ is the category of coherent models of the Heyting category $\mathcal{C}$. This is worked out in the author PhD thesis for the more general disjunction-free fragment.} which in particular means that the interpretation of $S$ in the presheaf topos will be precisely that of $\neg R$. Note that, if we add to the $(2^{\kappa})^+$-$Reg_{\neg}$ axiomatization above all instances of excluded middle for atomic formulas, we get an axiomatization of (an isomorphic copy of) the $\mu$-AEC, a fact which we will use in the following:

\begin{thm}\label{tameness0}
Let $\mathcal{K}$ be a $\mu$-AEC with directed colimits which is categorical in $\kappa$ and in $\lambda>2^{\kappa}$. Then $\mathcal{K}$ is $(2^{\kappa}, <\lambda)$-tame.
\end{thm}

\begin{proof}
Note first that we can take the model $M$ of size $\lambda$ as a monster model for $\mathcal{K}_{\geq (2^{\kappa})^+, < \lambda}$, since by the arguments of \cite{espindolas}, we have amalgamation there. Now Galois types in $M$ correspond to $\lambda^+$-geometric syntactic types, as shown in \cite{espindolas} (indeed, Galois types over $M_0$ correspond to syntactic types containing the complete formula that realizes the type of the tuple given by the underlying set of $M_0$). Thus, it is enough to show that a $\lambda$-coherent existential sentence of the form $\exists \mathbf{x} \phi(\mathbf{x}, \mathbf{d}, \mathbf{c})$, with constants $\mathbf{c}$ from the submodel $M_0$, the set of parameters of the type, where $\mathbf{d}$ is a finite tuple and where $\phi$ is a conjunction of atomic formulas, holds in $M$ if (and only if) every $(2^{\kappa})^+$-small approximation $\exists \mathbf{x}' \psi(\mathbf{x}', \mathbf{d}, \mathbf{c}')$ holds there. So suppose this latter condition holds. Let $N$ be a $(2^{\kappa})^+$-pure submodel containing $\mathbf{c}$ and $\mathbf{d}$, and consider the following theory: to the $(2^{\kappa})^+$-$Reg_{\bot}$ Morleyization of the sequents in $Reg_{\neg}$ logic that axiomatize $\mathcal{K}$, add the diagram of $N$, sequents expressing those negated existential sentences with constants from $N$ holding there, and sequents expressing that the $(2^{\kappa})^+$-small approximations $\exists \mathbf{x}' \psi(\mathbf{x}', \mathbf{d}, \mathbf{c}')$ hold. Clearly, every $(2^{\kappa})^+$-small subset has a model (the obvious expansion of the monster model) and so the whole theory has a $(2^{\kappa})^+$-$Reg_{\bot}$ model. This means that there is a $(2^{\kappa})^+$-pure embedding of $N \to S$ into a $(2^{\kappa})^+$-$Reg_{\bot}$ model of $\exists \mathbf{x} \phi(\mathbf{x}, \mathbf{c})$. By a similar proof to that of Grossberg conjecture in \cite{espindolas}, we can see that $M$ is injective with respect to $(2^{\kappa})^+$-pure embeddings in $\mathcal{K}^r_{\geq (2^{\kappa})^+, < \lambda}$, and thus we get that $\exists \mathbf{x} \phi(\mathbf{x}, \mathbf{d}, \mathbf{c})$ must hold there, as we wanted to prove.  
\end{proof}

  Assuming categoricity in a sufficiently large initial segment, we can also derive tameness:

\begin{thm}\label{tameness}
Let $\mathcal{K}$ be a $\mu$-AEC with directed colimits which is $\mu$-categorical for $\kappa \leq \mu < \beth_{\omega}(\kappa)$. Then $\mathcal{K}$ is $2^{\kappa}$-tame. In particular, categoricity in $[\kappa, \beth_{\omega}(\kappa))$ implies categoricity everywhere above $\kappa$.
\end{thm}

\begin{proof}
Note first that our whole analysis before Theorem \ref{tameness0} could be upgraded to the case in which we know that we have categoricity in $\mu$ for $\kappa \leq \mu < \chi := \beth_{\omega}(\kappa)$. In this case, it is possible to have an axiomatization in $\chi^+$-$Reg_{\neg}$ logic by adding a $\mu^+$-arity predicate $P_{\mu}$ for each $\mu<\chi$ and proceeding similarly to the axiomatization above. Let $\mathbb{T}^r$ be the $\chi^+$-$Reg_{\bot}$ Morleyization of the the following theory in the disjunction-free fragment: to the $\chi^+$-$Reg_{\neg}$ axiomatization of the models of size at least $\chi$, we add all instances of the axioms $\bigwedge_{i<\chi}\neg \neg \phi_i \vdash_{\mathbf{x}} \neg \neg \bigwedge_{i<\chi} \phi_i$, where $\phi_i$ are $<\chi$-$Reg_{\bot}$ formulas. Let also $\mathcal{K}_{\geq \chi, < \lambda}^{B, r}$ be its category of $\chi^+$-$Reg_{\bot}$ models of size at least $\chi$ and less than $\lambda$ with $\chi$-Boolean homomorphisms. It is enough to prove that the subtopos $\Sets^{\mathcal{K}_{\geq \chi, < \lambda}} \hookrightarrow \Sets^{\mathcal{K}_{\geq \chi, < \lambda}^{B, r}}$ is dense. As in the proof of Grossberg conjecture from \cite{espindolas}, we know similarly that this subtopos is obtained by adding the $\lambda$-topology generated by instances of excluded middle, so we just need to prove that the sequents $\bigwedge_{i<\lambda}\neg \neg \phi_i \vdash_{\mathbf{x}} \neg \neg \bigwedge_{i<\lambda} \phi_i$ hold in $\Sets^{\mathcal{K}_{\geq \chi, < \lambda}^{B, r}}$, for which it is in turn enough to prove that $\mathcal{K}_{\geq \chi, < \lambda}^{B, r}$ has directed bounds. Since it has bounds of chains of cofinality at least $\chi^+$, it is enough to consider chains of cofinality less than $\chi^+$. So let $\{M_i\}_{i<\alpha}$ be such a chain; we need to prove that there is a model containing the union of the diagrams in the chain. By $\chi^+$-$Reg_{\bot}$ compactness, we can assume without loss of generality that all models in the chain have size at most $\chi$, in which case it is enough to take a $\chi^+$-saturated model as a directed bound.\footnote{If the chain has cofinality bigger than $\omega$, at limit levels we use the existence of weakly initial models of the $\chi^+$-$Reg_{\bot}$ theory of the union of models below that level.} To show such a model exists, we will prove in the next paragraph that $\mathcal{K}_{\chi}^{B, r}$ has amalgamation at $\chi$. On the other hand, note that the double negation subtopos of $\Sets^{\mathcal{K}_{\chi}^{B, r}}$ satisfies all sequents in $\mathbb{T}^r$ due to axioms added to the axiomatization, and that it has a $\chi^+$-point since it $\chi^+$-classifies a $\chi^+$-$Reg_{\bot}$ theory: it is obtained as the quotient of $\mathbb{T}^r$ by all sequents expressing that morphisms $[M, -] \to [M', -]$ are epimorphisms. Thus, if $M$ is a $\chi^+$-point of the double negation subtopos, i.e., a $\chi^+$-saturated model of $\mathbb{T}^r$, then $M$ is clearly a directed bound for the given chain, as we wanted. 

  It remains to prove our claim. This can be seen through an argument using $3$-amalgamation of models of $\mathbb{T}^r_n := \mathbb{T}^r \cap \mathcal{L}_{\beth_{n}(\kappa)^+}$ of size $\beth_{n}(\kappa)^+<\chi$ and $\beth_{n}(\kappa)^+$-pure morphisms. First, $3$-amalgamation follows by taking as the amalgam at each level $\beth_{n}(\kappa)^+$ a weakly initial model of the pushout, in the doctrine of $\beth_{n}(\kappa)^+$-coherent categories, of the $\beth_{n}(\kappa)^+$ theories of the models that constitute the amalgamation diagram, which must be consistent since amalgamation holds. Since the $\beth_{n}(\kappa)^+$-theory of each pushout is axiomatized by $\beth_{n}(\kappa)^+$-$Reg_{\bot}$ axioms, they have weakly initial models $M_n$. Let us see that we can choose the homomorphisms between them to form a directed diagram whose colimit diagram $D$ is consistent with $\mathbb{T}^r$. Indeed, the diagram of the smallest weak initial model is consistent with $\mathbb{T}^r$ and so it has a model $M$; it suffices then to add to each $\beth_{n}(\kappa)^+$-theory of the pushouts also the diagrams of the weakly initial models of size $\beth_{n-1}(\kappa)^+$ that it contains: this provides a canonical homomorphism and we can use again weak initiality of these and successively embed them into $M$. Whence $M$ is the desired amalgam. This concludes the proof.
\end{proof}

\section{Shelah categoricity conjecture for AEC's}\label{scc}

We now get to the following:

\begin{thm}\label{aec}
(Shelah categoricity conjecture for AEC's). Let $\mathcal{K}$ be an AEC. If $\mathcal{K}$ is categorical in some $\lambda \geq \beth_{(2^{LS(\mathcal{K})})^+}$, then $\mathcal{K}$ is $\lambda'$-categorical for every $\lambda' \geq \beth_{(2^{LS(\mathcal{K})})^+}$. 
\end{thm}

\begin{proof}
Assume that the AEC is $\lambda$-categorical for some $\lambda \geq \beth_{(2^{LS(\mathcal{K})})^+}$. Consider the double negation subtopos of $\Sets^{\mathcal{K}_{\kappa}}$, corresponding to the dense topology, when $\kappa=LS(\mathcal{K})$. The first observation is that the dense topology is generated by those covers $\{M \to N_i\}$, where the $N_i$ have the same cardinality of $M$, such that any $M \to N$, for an arbitrarily large model $N$, can be amalgamated with one $M \to N_j$. Indeed, this choice of morphisms is enough to contain all covers $\{M \to N_i\}_{i \in I}$ where, for some partition $I = J \cup K$, $\{[N_j, -] \rightarrowtail [M, -]\}_{j \in J}$ corresponds to existential $\nu$-coherent sentences whose union is $\phi$ and $\{[N_k, -] \rightarrowtail [M, -]\}_{k \in K}$ corresponds to those whose union is $\neg \phi$. This guarantees that the sheaf topos is Boolean (whence, it is the double negation subtopos). Clearly the model of size $\lambda$ must be one of its points, so there is a type in the Boolean algebra of subterminal objects corresponding to the sentences which hold in the model of size $\lambda$. Since it is also Boolean, it is atomic, and thus the type is isolated by a subobject $S \rightarrowtail 1$. It follows that the slice over $S$ is a two-valued Boolean topos all of whose $\kappa^+$-points are $\mathcal{L}_{\infty, \kappa}$-elementarily equivalent; moreover, since there is a $\kappa^+$-saturated model of size $\kappa^+$ by stability (which holds from categoricity at $\lambda$), and any such model must be a $\kappa^+$-point of the two-valued Boolean subtopos, any two such $\kappa$-points are isomorphic. This means that the $\mathcal{L}_{\kappa^{++}, \kappa^+}$-theory $\kappa^+$-classified by the slice subtopos defines a $\kappa^+$-AEC categorical at $\kappa^+$ and $\lambda$. By Lemma \ref{dc} it must be closed under directed colimits and is, thus, an AEC. Indeed, we do not need to use amalgamation, since $\kappa^+$-closed models are precisely the $\kappa^+$-points of the double negation subtopos, and if they satisfy in addition the same sentences of the $\kappa^+$-saturated model they factor through the slice topos and are thus $\kappa^+$-saturated. As a consequence of Theorems \ref{tameness0} and \ref{tameness}, it is $\chi$-tame for $\chi = 2^{\kappa}$ and thus it is categorical in a proper class of cardinals. By the same arguments from \cite{espindolas}, we conclude that we have categoricity everywhere above $\kappa^+$. By Morley's omitting types theorem for AEC's\footnote{This is an application of Morley's method to the AEC to prove the analogous version of Proposition 3.3 in \cite{shma}.}, we conclude that the every model of the AEC is $\kappa^+$-saturated above some $\chi<\beth_{(2^{LS(\mathcal{K})})^+}$. Our theorem follows.
\end{proof}

  The threshold cardinal is best possible, since by one of Shelah's example (see e.g. Fact 9.8 in \cite{sv3}) for any $\alpha<(2^{LS(\mathcal{K})})^+$ there are AEC's with models of size up to $\beth_{\alpha}(LS(\mathcal{K}))$ but no bigger models. This allows to construct AEC's that are eventually categorical but non categorical above $\beth_{\alpha}(LS(\mathcal{K}))$. The threshold in the special case of $\mathcal{L}_{\kappa^+, \omega}$ has been conjectured by Shelah and appears also in his work from the nineties.

\section{An infinitary extension of Morley's categoricity theorem}\label{im}

We can now proceed with the:

\textit{Proof of Theorem \ref{morley}}: Assume we have an accessible category with directed colimits $\mathcal{K}$ which is categorical for some $\lambda \geq \beth_{(2^{\kappa})^+}$ in $S$. By Theorem 4.5 in \cite{br}, $\mathcal{K}$ is a reflective subcategory of a finitely accessible category, $\mathcal{K}'$, which must be equivalent to the category of models of a $\omega$-geometric theory; without loss of generality we can restrict our accessible category to monomorphisms and so we can assume the morphisms of $\mathcal{K}'$ are embeddings. Both categories are axiomatizable in $\mu$-coherent logic for $\mu=LS(\mathcal{K})^+$. Let $\lambda$ be the first cardinal in $S$ bigger than $LS(\mathcal{K})$; it follows that the inclusion $i: \mathcal{K}_{\geq \mu, < \lambda} \hookrightarrow \mathcal{K}_{\geq \mu, < \lambda}'$ has a left adjoint $\pi$. Note that $\pi$ is computed by considering the expression of $M$ as a directed colimit of finitely presentable models and setting $\pi(M)$ to be the directed colimit in $\mathcal{K}$ of the same diagram. It further follows that $i^*: \Sets^{\mathcal{K}_{\geq \mu, < \lambda}'} \to \Sets^{\mathcal{K}_{\geq \mu, < \lambda}}$ is the inverse image of a geometric morphism which is left adjoint to $\pi^*$.\footnote{Incidentally, that $i^*$ is right adjoint to $\pi^*$ is also true since it follows from the adjunction $\pi \dashv i$.}. To see this, we need to check that there is an isomorphism $[i^*(F), G] \cong [F, \pi^*(G)]$ natural in $F$ and $G$, and since $i^*$ preserves colimits, it is enough to take $F \cong [A, -]$ a representable. Then an element of $[i^*(F), G]$ is a compatible family of functions $[A, i(N)] \to G(N)$ indexed by $\mathcal{K}$-models $N$, or, what amounts to the same thing, a compatible family $[A, \pi(M)] \to G(\pi(M))$ indexed by $\mathcal{K}'$-models $M$. But by naturality this latter family determines completely a compatible family $[A, M] \to G(\pi(M))$ (given that $\pi(\pi(M))=\pi(M)$), which is the same as an element in $[F, \pi^*(G)]$. Note as well that $\pi^*$ must be an embedding which is clearly dense, and so both presheaf toposes will have the same $\lambda$-saturated models.

  Consider now the $\kappa^+$-AEC $\mathcal{K}''$ consisting of $\kappa^+$-saturated models. By Lemma \ref{dc}, we can see that $\mathcal{K}''$ inherits the concrete directed colimits of $\mathcal{K}'$. Also, since Galois types over submodels of size $\mu$ correspond to $\mu^+$-geometric types containing the complete formula satisfied by the underlying set of the submodel, we can use the same proof idea of Proposition 3.5 from \cite{shma} to derive stability from categoricity of $\mathcal{K}$ in $\lambda$, for which it is enough to note that the usual Ehrenfeucht-Mostowski model $M$ with $\lambda$ indiscernibles, which is a directed colimit (in $\mathcal{K}$) of those with finitely many indiscernibles (which are as well models of $\mathcal{K}$) omits a type over some submodel if and only if $\pi(M)$ omits such type. Indeed, the directed colimit (in $\mathcal{K}''$) must embed into each one of the models which is a Skolem hull of the directed union of models, and thus a realization of the type there would entail a realization of the type in all possible Skolem hulls, whence by the transfinite transitivity rule up to double negation the type would be already realized in the directed union (i.e. in $M$).\footnote{Since our theories and types are $\lambda$-coherent, we can eliminate the double negation by conservativity of the classical fragment over the coherent fragment} Finally, by Theorem \ref{aec} applied to the AEC $\mathcal{K}''$ we get thus categoricity in every $\lambda \geq \kappa^+$, and since by Morley's omitting types theorem applied to $\mathcal{K}'$ we have that all models of size at least $\nu$ for some $\nu < \beth_{(2^{\kappa})^+}$ belong to $\mathcal{K}''$, we get in particular that $\mathcal{K}$ is categorical in those cardinals in $S$ above $\beth_{(2^{\kappa})^+}$, as we wanted.
  
  Finally, we consider the case in which the sentence defining the category of models is a compact sentence, and we consider those embeddings that are first-order elementary. This is a $\mu$-AEC with directed colimits, amalgamation and no maximal models. The corresponding finitely accessible category $\mathcal{K}'$ will also have, then, amalgamation and no maximal models. As we have proven above, $\mathcal{K}'$ will be eventually categorical. Since by the considerations in the proof of Lemma \ref{dc} we have a stable surjection $\Sets^{M/\mathcal{K}_{\geq \kappa, < \lambda}^B} \twoheadrightarrow \Sets^{M/\mathcal{K}_{\geq \kappa, < \lambda}}$, each model $M$ of $\mathcal{K}'$ has a proper $\omega$-pure (in fact, $\omega$-Boolean) extension, as otherwise $\Sets^{M/\mathcal{K}_{\geq \kappa, < \lambda}^B}$ would be two-valued and Boolean, forcing $\Sets^{M/\mathcal{K}_{\geq \kappa, < \lambda}}$ to be two-valued and Boolean, which is impossible since $M$ is not maximal. It follows that each model has a $\omega$-pure morphism into the model of size $\lambda$, and thus every such model is $\omega$-saturated; in particular, $\omega$-coherent formulas are either $0$ or $1$. By the arguments in section $9$ of \cite{espindolas}, any sentence of the form $\forall \mathbf{x} (\vartheta \to \eta)$, where $\vartheta$ and $\eta$ are $\omega$-coherent, which is valid in the $LS(\mathcal{K})^+$-saturated model is provable in the $LS(\mathcal{K})^+$-classifying topos of models of size at least $LS(\mathcal{K})^+$ (note that $\forall \mathbf{x} (\vartheta \to \eta)$ will also be equivalent to a $\omega$-coherent formula, as can be seen  by compactness arguments using the definition of universal quantification in the syntactic category). This readily implies, as can be seen by the proof in \cite{keisler} for conjunctive formulas in saturated models, that any $\omega_1$-coherent formula is equivalent to the conjunction of their approximations (in particular, $\omega_1$-coherent sentences are either $0$ or $1$). This allows thus to prove that now sentences of the form $\forall \mathbf{x} (\vartheta \to \eta)$, where $\vartheta$ and $\eta$ are $\omega_1$-coherent, which are valid in the $LS(\mathcal{K})^+$-saturated model is provable in the $LS(\mathcal{K})^+$-classifying topos. Continuing with this process, we finally get that all $LS(\mathcal{K})^+$-coherent sentences are either $0$ or $1$, which is enough to get categoricity at $LS(\mathcal{K})^+$ and hence everywhere above $LS(\mathcal{K})^+$. This finishes the proof.

\section{Classification of categoricity spectra}  
  
We end with:

\begin{thm}\label{cl}
Let $\mathcal{K}$ be a large $\kappa$-accessible category with directed colimits. Assume the Singular Cardinal Hypothesis $SCH$ (only if the restriction to monomorphisms is not an AEC). Then the categoricity spectrum $\mathcal{C}at(\mathcal{K})=\{\lambda\geq \kappa: \mathcal{K} \text{ is $\lambda$-categorical}\}$ is one of the following:

\begin{enumerate}
\item $\mathcal{C}at(\mathcal{K})=\emptyset$.
\item $\mathcal{C}at(\mathcal{K})=[\alpha, \beta]$ for some $\alpha, \beta \in [\kappa, \beth_{\omega}(\kappa))$.
\item $\mathcal{C}at(\mathcal{K})=[\chi, \infty)$ for some $\chi \in [\kappa, \beth_{(2^{\kappa})^+})$.
\end{enumerate}

\end{thm}  

\begin{proof}
Our proof idea shares the same guidelines as the one in \cite{sv3}, except that the amalgamation hypothesis is not used and $WGCH$ is replaced with $SCH$ or, in AEC's, eliminated completely (all this at the price of replacing $\kappa^{+\omega}$ with $\beth_{\omega}(\kappa)$). If case $1$ does not occur, suppose first there is some categoricity cardinal $\lambda \geq \beth_{\omega}(\kappa)$, and proceed as in the proof of Theorem \ref{aec} to define a $\kappa^+$-AEC with directed colimits categorical in $\kappa^+$ and $\lambda$. By Theorems \ref{tameness0} and \ref{tameness}, it must be $2^{\kappa^+}$-tame and so must be the original category of models of $\phi$, which allows us to conclude that case $3$ occurs. If categoricity occurs only below $\beth_{\omega}(\kappa)$, it must be a segment as proven in \cite{espindolas}. Indeed, we do not need amalgamation by using the same argument as in the proof of Theorem \ref{aec}. On the other hand, the assumption of maximal models in the proof of Theorem 9.1 in \cite{espindolas} can be dropped by building the tree of theories so that those theories $\Gamma$ that admit a maximal model are not extended to $\Gamma'$ by adding a new constant and so they skip this modification of the construction; in the end, in the tree of theories, branches reach to the level $\lambda$ or they reach a node in which the theory is classical (corresponding to a maximal model). In either case, instances of excluded middle can be shown to hold at the topos $\Sets[\mathbb{T}_{\kappa^+}]_{\kappa^+}$, making it Boolean. To sum up, categoricity cannot alternate, which leaves us in the case $2$ and thus finishes the proof.
\end{proof} 

Examples culled from the literature on AEC's show that each of the three cases in the classification can indeed occur (see e.g. the examples of \cite{sv3}). The non-trivial case is $2$, which occurs in the Shelah-Villaveces example in \cite{shvi}.

\bibliographystyle{amsalpha}

\renewcommand{\bibname}{References} % changes the header from Bibliography to References

\bibliography{references}

%\begin{thebibliography}{widest entry}

%\end{thebibliography}

\end{document}